\input amstex
\documentstyle{amsppt}
\magnification1200
\pagewidth{6.5 true in}
\pageheight{9.25 true in}
\NoBlackBoxes

\topmatter
\title 
Partial sums of the M{\" o}bius function
\endtitle
\author 
 K. Soundararajan
\endauthor 
 
\address
{Department of Mathematics, 450 Serra Mall, Bldg. 380, Stanford University, 
Stanford, CA 94305-2125, USA}
\endaddress
\email 
ksound{\@}stanford.edu
\endemail
\thanks   The author is partially supported by the National Science Foundation (DMS 0500711) 
and the American Institute of Mathematics (AIM). 
\endthanks
 \endtopmatter

\document

\head 1. Introduction \endhead 
  
\noindent This paper is concerned with bounding
$$
M(x) = \sum_{n\le x} \mu(n).
$$ 
J.E. Littlewood [6] proved that if the Riemann Hypothesis (RH) is true then, for any 
fixed $\epsilon>0$, $1/\zeta(1/2+\epsilon+it) \ll |t|^{\epsilon}$.  
It follows by Perron's formula that 
$$
M(x) \ll x^{\frac 12+\epsilon}. \tag{1}
$$
Conversely, the estimate $M(x) \ll x^{\frac 12+\epsilon}$ implies, by partial summation, the 
convergence of the series $\sum_{n=1}^{\infty} \mu(n)n^{-s} =1/\zeta(s)$ for any $\sigma>1/2$, 
and therefore RH.  Subsequently, E. Landau [5] showed that, assuming RH, (1) is valid with 
$\epsilon \ll \log \log \log x/\log \log x$, and E.C. Titchmarsh [13] improved this to $\epsilon \ll 
1/\log \log x$.  H. Maier and H.L. Montgomery [7] obtained a substantial improvement over 
these results, and established that 
$$
M(x) \ll x^{\frac 12} \exp\Big( C(\log x)^{\frac {39}{61}}\Big). \tag{2}
$$
They comment that the limit of their method would be an exponent in (2) slightly 
smaller than $39/61$.  In this paper, we introduce some new ideas which 
permit the following better result. 

\proclaim{Theorem 1}  Assume RH.  For large $x$ we have
$$
M(x)\ll  \sqrt{x} \exp( (\log x)^{\frac 12}(\log \log x)^{14}).
$$
\endproclaim 
     
The main ingredient in our proof is a result on the frequency with which 
ordinates of the zeros of $\zeta(s)$ can cluster in short intervals, which may be of 
independent interest. 
Let $N(T)$ denote the number of zeros $\rho=\beta+i\gamma$ 
of $\zeta(s)$ with ordinate $\gamma$ lying in $[0,T]$.  Recall that 
$$
N(T) = \frac{T}{2\pi} \log \frac{T}{2\pi e} + \frac 78 +S(T) +O\Big(\frac 1T\Big), 
$$
where $\pi S(T) =\text{arg} \zeta(\tfrac 12+iT)$ and the argument 
is obtained by continuous variation from $2$ (where the argument is zero) 
to $2+iT$ to $\tfrac 12+iT$.  It is easy to show that $S(T)\ll \log T$, and 
on RH Littlewood established that $S(T) \ll \log T/\log\log T$.  Recently D. Goldston and 
S. Gonek [3] put Littlewood's bound into the elegant form $|S(T)| \le (\tfrac 12+o(1)) 
\log T/\log \log T$.  Building on their work, we quantify here the frequency of large values of 
$S(t+h) -S(t-h)$; equivalently, the frequency with which the interval 
$[t-h,t+h]$ contains an unusual number of ordinates of zeros of $\zeta(s)$.

\proclaim{Theorem 2} Assume RH.  Let $T$ be large, and 
let $0\le h \le \sqrt{T}$, and $(\log \log T)^2 \le V\le \log T/\log \log T$ be given.  
The number of well-spaced points $T\le t_1 <t_2 <\ldots < t_R \le 2T$ with 
$t_{j+1} -t_j \ge 1$ and such that 
$$
\Big| N(t_j+h)-N(t_j-h) - \frac{h}{\pi} \log \frac{t_j}{2\pi } \Big|  >V 
$$
satisfies the bound 
$$
R \ll T \exp\Big(-V\log \frac{V}{\log \log T} +3 V\log \log V\Big).
$$
\endproclaim  

In [10, 11] A. Selberg established unconditionally that $\pi S(t)$ has a Gaussian distribution with mean 
$0$ and variance $\frac 12 \log \log T$.   This suggests a better bound for $R$ than 
that furnished by Theorem 2.  Namely, perhaps the bound 
$R\ll T\exp(-C V^2/\log \log T)$ holds for some absolute positive constant $C$, uniformly 
in $V$.   This is in keeping with the recent conjecture of D.W. Farmer, Gonek and C.P. Hughes 
[2] that $S(t) \ll \sqrt{\log T \log \log T}$.  
By adapting the ideas in [12] it would be possible to establish the conjectured 
bound for $R$ (assuming RH) in the range $V\ll (\log\log T) \log \log \log T$. A more detailed analysis 
of such results is the focus of my ongoing joint work with Chris Hughes and 
Nathan Ng.

Using Theorem 2 we shall establish an estimate
for the frequency with which small values of $|\zeta(s)|$ are 
attained.   The main result of my paper [12] deals with 
corresponding estimates for the frequency with which large 
values of $|\zeta(\tfrac 12+it)|$ are attained.  To state our results 
conveniently we require a definition. 

\proclaim{Definition 3}  Let $T$ be large and let $(\log \log T)^2 \le 
V\le \log T/\log \log T$ be given.
We say that a point $t\in [T,2T]$ is {\rm $V$-typical} if the following three 
conditions hold; if one of these criteria fails, we say that the 
point is {\rm $V$-atypical}.  

\noindent (i).  Let $x=T^{1/V}$.  For all $\sigma \ge \tfrac 12$ we have 
$$
\Big| \sum_{n\le x} \frac{\Lambda(n)}{n^{\sigma+it} \log n} \frac{\log (x/n)}{\log x} 
\Big| \le 2V. 
$$

\noindent (ii).  Every sub-interval of $(t-1,t+1)$ of length $2\pi V/\log T$ 
contains at most $3V$ ordinates of zeros of $\zeta(s)$. 

\noindent (iii).  Every sub-interval of $(t-1,t+1)$ of length $2\pi V/((\log V)\log T)$ 
contains at most $V$ ordinates of zeros of $\zeta(s)$.  
\endproclaim

\proclaim{Proposition 4} Assume RH.  Let $T$ be large.  Any point 
$t \in [T,2T]$ is $V$-typical provided 
$V\ge (\tfrac 12+o(1)) \log T/\log \log T$.   Given $(\log \log T)^2\le V\le \log T/\log \log T$, 
the number of well-spaced $V$-atypical points $T\le t_1 \le \ldots \le t_R \le 2T$ with 
$t_{j+1}-t_j \ge 1$  satisfies 
$$
R\ll T \exp\Big(- V\log \frac{V}{\log \log T} + 4V\log \log V\Big). 
$$ 
\endproclaim 

\proclaim{Proposition 5}  Assume RH.  Let $T$ be large, and suppose
$t \in [T,2T]$ is $V$-typical for some $(\log \log T)^2 \le V\le \log T/\log \log T$.  
  Put 
$\sigma_0= \tfrac 12+ \tfrac V{\log T}$.   For $2\ge \sigma \ge \sigma_0$ 
we have 
$$
\log |\zeta(\sigma+it)| \ge - V\log\log V, \tag{3}
$$ 
and for $\tfrac 12 < \sigma \le \sigma_0$ we have 
$$
\log |\zeta(\sigma+it)| \ge - V\log \Big(\frac{\sigma_0-1/2}{\sigma-1/2} \Big) -  8V\log \log V. \tag{4}
$$
\endproclaim 

We will describe in \S 5 below how our main result, Theorem 1, 
follows from Propositions 4 and 5, and a careful application of Perron's formula.  Just as 
we expect that the true bound for $R$ in Theorem 2 should be much smaller, we may 
expect a corresponding improvement of Proposition 4.  Perhaps 
the better bound $R\ll T\exp(-CV^2/\log \log T)$ holds, for some positive constant $C$.  
If such were the case, then our method would yield 
$M(x)\ll x^{\frac 12} \exp(C(\log \log x)^{3})$ for some positive constant $C$.   Even this is 
far from the conjectured maximal order of magnitude for $M(x)$: Gonek (unpublished, but see 
N. Ng [8]) has conjectured that 
$$
\infty >\limsup_{x\to \infty} \frac{M(x)}{\sqrt{x} (\log \log \log x)^{\frac 54}}>0 > \liminf_{x\to \infty} 
\frac{M(x)}{\sqrt{x}(\log \log \log x)^{\frac 54}}>-\infty.
$$

{\bf Acknowledgments.}  I am grateful to Professors Maier and Montgomery 
for making available their preprint [7] which motivated the present paper.   I am 
also grateful to Gergely Harcos for a query which led to a clarification of the proof.
Part of the paper was written while I visited the University of Bristol.  I 
am most grateful to them for their kind hospitality.

\head 2.  Preliminary Lemmas \endhead 

\noindent We collect here three familiar results that 
we shall need below.  These are Selberg's construction of 
good approximations to characteristic functions of 
intervals, the explicit formula connecting 
primes and zeros, and a version of the large sieve. 

\proclaim{Lemma 6} Let $h >0$ and $\Delta >0$ be given.  Let $\chi_{[-h,h]}$ denote the 
characteristic function of the interval $[-h,h]$.  There 
exist even analytic functions $F_{-}$, and $F_+$ satisfying the following properties. 

(i) $F_{-}(u) \le \chi_{[-h,h]}(u) \le F_+(u)$ for real $u$.

(ii) We have 
$$
\int_{-\infty}^{\infty} |F_{\pm}(u)-\chi_{[-h,h]}(u)| du \le 1/\Delta.
$$

(iii) ${\hat F}_{\pm}(x)=0$ for $|x| \ge \Delta$ where ${\hat F}_{\pm}(x)= \int_{-\infty}^{\infty} 
F_{\pm}(u) e^{-2\pi ixu} du$ denotes the Fourier transform. Also, 
$$
{\hat F}_{\pm}(x) = \frac{\sin (2\pi hx)}{\pi x} + O\Big( \frac 1{\Delta}\Big). 
$$

(iv) If $z=x+iy$ is a complex number with $|z| \ge 2h$ then
 $$
 |F\pm(z)| \ll \frac{e^{2\pi \Delta |y|} }{\Delta |z|^2}.
 $$
\endproclaim 

\demo{Proof} Such functions were constructed by Selberg (see [9]), using Beurling's 
approximation to the signum function.  We give a brief description; for a detailed discussion see 
J.D. Vaaler [15].   Set $K(z)=(\sin \pi z)^2/(\pi z)^2$ and 
$$
H(z) = \Big( \frac{\sin \pi z}{\pi } \Big)^2 \Big( \sum_{n=-\infty}^{\infty} \frac{\text{sgn}(n)}{(z-n)^2} 
+\frac{2}{z}\Big), 
$$ 
where $\text{sgn}(x)$ is the sign function taking values $1$ for positive $x$, $-1$ for negative 
$x$, and $0$ for $x=0$.   Beurling showed that $H(x) -K(x) \le \text{sgn}(x) \le H(x)+K(x)$, 
and that 
$$
\int_{-\infty}^{\infty} |H(x) \pm K(x) - \text{sgn}(x)| dx =1.
$$
The desired functions $F_\pm$ are given by 
$$
F_{\pm}(z) = \frac 12 \Big( H(\Delta(x+h)) \pm K(\Delta(x+h)) + H(\Delta(h-x)) \pm K(\Delta(h-x))\Big).
$$
Properties (i)-(iii) are well-known, and it is not difficult to check the bound in (iv).  
\enddemo 

\proclaim{Lemma 7}  Let $h(s)$ be analytic in the strip $|\text{Im}(s)| \le \tfrac 12+\epsilon$ 
for some $\epsilon >0$, taking real values on the real line,  and satisfying 
$|h(s)| \ll (1+|s|)^{-1-\delta}$ for some $\delta >0$.  Then, with $\rho =\tfrac 12+i\gamma$ 
denoting the non-trivial zeros of $\zeta(s)$,  
$$
\align
\sum_{\rho} h(\gamma) &= h\Big( \frac1{2i}\Big) + h\Big(-\frac 1{2i} \Big)   
+ \frac 1{2\pi} \int_{-\infty}^{\infty} h(u) \Big( \text{Re} \frac{\Gamma^{\prime}}{\Gamma} 
\Big(\frac 14+ \frac{ iu}{2} \Big) -\log \pi \Big) du\\
& - \frac{1}{2\pi } 
\sum_{n=2}^{\infty} \frac{\Lambda(n)}{\sqrt{n}} \Big( {\hat h}\Big( \frac{\log n}{2\pi} \Big) 
+ {\hat h} \Big( -\frac{\log n}{2\pi}\Big) \Big). 
\\
\endalign
$$
\endproclaim 

\demo{Proof}  This is the explicit formula; see for example Lemma 1 in [3], or Chapter 
5 of H. Iwaniec and E. Kowalski [4].
\enddemo  
 
\proclaim{Lemma 8} Let $A(s) = \sum_{p\le N} a(p)p^{-s}$ be a Dirichlet 
polynomial.  Let $T$ be large and suppose $s_r= \sigma_r+it_r$ ($r=1$, $\ldots$, $R$) 
be points with $T <t_1 < t_2 < \ldots <t_R \le 2T$ and $t_{r+1}-t_r \ge 1$, 
and $\sigma_r \ge \alpha$.  For any $k$ with $N^k \le T$ we have 
$$
\sum_{r=1}^{R} |A(s_r)|^{2k} \ll T(\log T)^2 k! \Big( \sum_{p\le N} |a(p)|^2 p^{-2\alpha}\Big)^{k}.
$$
\endproclaim 

\demo{Proof}  This large sieve type inequality may be found as Lemma 5 in Maier and 
Montgomery [7].  
\enddemo 

\head 3. Proof of Theorem 2\endhead 

\noindent We use Lemma 6 to approximate the characteristic function of 
$[-h,h]$, taking there $\Delta = (1+\eta)(\log T)/(2\pi V)$ with $\eta=1/\log V$.  Let 
$F_\pm$ denote the functions produced in Lemma 6.   We now appeal to the 
explicit formula, Lemma 7, taking $h(s)= F_{\pm}(s-t)$ where $T\le t\le 2T$.  
Observe that ${\hat h}(x)= {\hat F_{\pm}}(x) e^{-2\pi i xt}$.  Therefore, 
the explicit formula gives
$$
\align
\sum_{\rho} F_\pm(\gamma-t) &= F_{\pm}\Big(\frac{1}{2i}-t\Big) + F_{\pm} \Big(-\frac{1}{2i}-t\Big) 
- \frac{1}{\pi}\text{Re } \sum_{n=2}^{\infty} \frac{\Lambda(n)}{n^{\tfrac 12+ it}} 
{\hat F}_{\pm}\Big(\frac{\log n}{2\pi}\Big)\\
&+ \frac{1}{2\pi} \int_{-\infty}^{\infty} F_{\pm}(u) \Big(\text{Re } \frac{\Gamma^{\prime}}{\Gamma} 
\Big(\frac 14+ i \frac{t+u}{2} \Big ) -\log \pi \Big)  du. \tag{5}
\\
\endalign
$$ 

Using Stirling's formula we may readily check that for $0< h\le \sqrt{T}$ 
(or see equation (13) of [3]) 
$$
\frac{1}{2\pi}\int_{-\infty}^{\infty} F_{\pm}(u) \Big(\text{Re } \frac{\Gamma^{\prime}}{\Gamma} 
\Big(\frac 14+ i \frac{t+u}{2} \Big ) -\log \pi \Big)  du = 
\frac{1}{2\pi} \log \frac{t}{2\pi} {\hat F}_{\pm}(0) + O(1).
$$
Note that with the $+$ choice of sign the LHS of (5) is at least $N(t+h)-N(t-h)$, 
while with the $-$ choice of sign it is at most $N(t+h)-N(t-h)$.  Moreover 
${\hat F}_{+}(0) \le 2h +1/\Delta$, and ${\hat F}_{-}(0) \ge 2h -1/\Delta$.  
These observations lead to 
$$
\align
\Big| N(t+h) - N(t-h) -\frac{h}{\pi} \log \frac{t}{2\pi}\Big| 
&\le \frac{1}{\Delta}\frac{\log T}{2\pi} + \max_{\pm} 
\Big( \Big| F_{\pm}\Big(\frac{1}{2i}-t\Big) \Big| +\Big| 
F_{\pm} \Big(\frac{-1}{2i}-t\Big)\Big|\\
& + \frac{1}{\pi}\Big|\sum_{n=2}^{\infty} 
\frac{\Lambda(n)}{n^{\frac 12+it}} {\hat F}_{\pm}\Big(\frac{\log n}{2\pi}\Big)\Big|\Big)+O(1). \tag{6}\\
\endalign
$$

Now $\log T/(2\pi \Delta) =V (1-\eta +O(\eta^2))$ and by part (iv) of Lemma 6 the 
contribution of $F_{\pm}(\pm 1/2i -t)$ terms is $\ll T^{-1}$.  Therefore if 
the LHS of (6) exceeds $V$ then we must have 
$$
\max_{\pm} \Big| \sum_{n=2}^{\infty} \frac{\Lambda(n)}{n^{\frac 12+it}} {\hat F}_{\pm} \Big( \frac{\log n}{2\pi} 
\Big) \Big| \ge  2\eta V.
$$
Since ${\hat F}_{\pm}(x)=0$ for $|x|\ge \Delta$ the sums above may be restricted to 
$n\le \exp(2\pi \Delta) = T^{(1+\eta)/V}$.  Moreover, the contribution of prime 
cubes and higher powers is $O(1)$.  Thus we have either 
$$
\max_{\pm} \Big| \sum_{p\le T^{(1+\eta)/V}} \frac{\log p}{p^{\frac 12+it}} 
{\hat F}_{\pm}\Big(\frac{\log p}{2\pi}\Big)\Big| \ge \eta V, 
\ \ \text{or} 
\ \ 
\max_{\pm} \Big| \sum_{p\le T^{(1+\eta)/2V} }
\frac{\log p}{p^{1+2it}} {\hat F}_{\pm}\Big(\frac{\log p}{\pi}\Big) \Big| \ge \eta V.
$$

We conclude, for our sequence of $R$ well-spaced points $t_j$, and any positive 
integer $k$ that 
$$
R (\eta V)^{2k} \le \sum_{j=1}^{R} \sum_{\pm} \Big(\Big| \sum_{p\le T^{(1+\eta)/V}} \frac{\log p}{p^{\frac 12+it}} 
{\hat F}_{\pm}\Big(\frac{\log p}{2\pi}\Big)\Big|^{2k} 
+ \Big|\sum_{p\le T^{(1+\eta)/2V} }
\frac{\log p}{p^{1+2it}} {\hat F}_{\pm}\Big(\frac{\log p}{\pi}\Big) \Big|^{2k}\Big).
$$
Suppose that $k \le V/(1+\eta)$, so that Lemma 8 applies.  In that case we 
obtain that 
$$
\align
R (\eta V)^{2k} \ll T (\log T)^2 k^k \Big(& \Big( \sum_{p\le T^{(1+\eta)/V}} 
\frac{\log^2 p}{p} \Big| {\hat F}_{\pm} \Big( \frac{\log p}{2\pi}\Big)\Big|^2 \Big)^k 
\\
&+ \Big( \sum_{p\le T^{(1+\eta)/2V}} \frac{\log^2 p}{p^2} \Big|{\hat F}_{\pm} \Big( \frac{\log p}{\pi}
\Big)\Big|^2 \Big)^k \Big).
\\
\endalign
$$
Using (iii) of Lemma 6 we conclude that the above is 
$$
\ll T (\log T)^2 (C k \log \log T)^k, 
$$ 
for some positive constant $C$.  Hence 
$$
R \ll T (\log T)^2 \Big ( \frac{C k \log \log T}{\eta^2 V^2}\Big)^k, 
$$
and the Theorem follows upon recalling that $\eta=1/\log V$, and taking the largest permissible 
value for $k$, namely $\lfloor V/(1+\eta) \rfloor$.

\head 4.  Lower bounds for $|\zeta(s)|$: Proof of Propositions 4 and 5\endhead 

\demo{Proof of Proposition 4}   If $V\ge (\tfrac 12+ \epsilon) \log T/\log \log T$ 
then $x= T^{1/V} \le (\log T)^{2-\epsilon}$ so that criterion (i) of Definition 3 is met.  
Moreover, Goldston and Gonek's estimate (see Theorem 1 of [3]) that 
for large $t$ and $0<h\le \sqrt{t}$ one has $|N(t+h)-N(t)-\frac{h}{2\pi} \log \frac{t}{2\pi}| \le 
(\tfrac 12+o(1))\log t/\log \log t$, readily shows that 
criteria (ii) and (iii) are also met.  Therefore $t$ is $V$-typical 
for $V\ge (\tfrac 12+o(1)) \log T/\log \log T$.  

We now obtain the bound for the number $R$ of well-spaced $V$-atypical 
points.   If a point is $V$-atypical then one of the criteria (i), (ii), or (iii) must 
be violated.   Appealing to Lemma 8 we may show 
(arguing exactly as in our proof of Theorem 2 above) that the 
number of well-spaced points for which condition (i) fails is 
$\ll T(\log T)^2 \exp(-(2+o(1))V \log (V/\log \log T))$.  Theorem 2 
shows that the number of well-spaced points for which (ii) fails is 
$\ll T \exp(-(2+o(1)) V\log (V/\log \log T))$ as well.  Theorem 2 also 
shows that the number of well-spaced points for which condition (iii) fails is 
$$
\ll T\exp\Big(- V\log \frac{V}{\log \log T} + 4V\log \log V\Big).
$$
Hence the bound for $R$ claimed in Proposition 4 follows.

\enddemo

\demo{Proof of Proposition 5}  Suppose that $t$ is $V$-typical, so that conditions (i)-(iii) of Definition 3 hold.  
We must now establish the estimates (3) and (4).  For $s=\sigma +it$ we write 
$$
F(s) = \sum_{\rho} \text{Re }\frac{1}{s-\rho} = \sum_{\rho} \frac{(\sigma-1/2)}{(\sigma-1/2)^2+(t-\gamma)^2}.
$$
By Stirling's formula and Hadamard factorization we have (see (2.12.7) of Titchmarsh [14], or Chapter 
12 of H. Davenport [1]) 
$$
\text{Re }\frac{\zeta^{\prime}}{\zeta}(s) = F(s) - \frac{1}{2} \log T +O(1). \tag{7}
$$ 

\proclaim{Lemma 9}  Let $T\le t\le 2T$ be $V$-typical.  
For $\frac 12 <\sigma \le \sigma_0 =\frac 12+ \frac{V}{\log T}$, 
we have 
$$
\log |\zeta(\sigma+it)| \ge \log |\zeta(\sigma_0+it)| - V\log \frac{(\sigma_0-\frac 12)}{(\sigma-\frac 12)} 
-7V\log \log V.  
$$
\endproclaim
\demo{Proof} Using (7) we 
see that 
$$
\align
\log |\zeta(\sigma_0+it)| - \log |\zeta(\sigma+it)| &= \int_{\sigma}^{\sigma_0} 
\text{Re }\frac{\zeta^{\prime}}{\zeta}(u+it) du \le \int_{\sigma}^{\sigma_0} F(u+it) du 
\\
&= \frac 12 \sum_{\gamma} \log \frac{(\sigma_0-\frac 12)^2+(t-\gamma)^2}{(\sigma-\frac 12)^2+(t-\gamma)^2}. \tag{8}
\\
\endalign
$$
We split the sum over $\gamma$ into various intervals.  First we have 
the range where $|t-\gamma|$ is below $2\pi V/((\log V)\log T)$.  Second we have the intervals  
$2\pi (n+1/\log V)V/\log T \le |t-\gamma| \le 2\pi (n+1+1/\log V)V/\log T$ 
for $0\le n\le N=[(\log T)/(4\pi V)]$.  Finally there is the range $|t-\gamma| > 2\pi (N+1+1/\log V)V/\log T$.  
Using condition (iii) of Definition 3, we see that 
the first range contributes to (8) an amount $\le V \log ((\sigma_0-\frac 12)/(\sigma -\frac 12))$.  
In the second range we use condition (ii) of Definition 3, and conclude that the contribution 
of such terms to (8) is 
$$
\le {3V}\sum_{n=0}^{N} \log \frac{1+(n+1/\log V)^2}{(n+1/\log V)^2} 
\le 6V \log \log V + 10V.
$$
Splitting into intervals of length $1$, we see easily that 
the final range contributes 
$$
\le \frac{1}{2} \sum_{|t-\gamma|> 1/2} \frac{(\sigma_0-\frac 12)^2}{|t-\gamma|^2} =o(V).
$$
Putting everything together we obtain the Lemma. 
\enddemo

From Lemma 9, estimate (4) would follow once (3) is 
established.  In other words, we now need to deal with $\sigma \ge \sigma_0$.  For this we need the 
following Lemma. 

\proclaim{Lemma 10}  Let $t$ be large and let $T\le t\le 2T$.  
Uniformly for $\frac 12 < \sigma \le 2$, and $2\le x\le T$ we have 
$$
\log |\zeta(\sigma +it)| \ge \text{Re }\sum_{n\le x} 
\frac{\Lambda(n)}{n^{\sigma+it} \log n} \frac{\log (x/n)}{\log x} 
- \Big(1+ \frac{x^{\frac 12-\sigma}}{(\sigma-\frac 12) \log x}\Big) \frac{F(\sigma+it)}{\log x} +O(1). 
$$ 
\endproclaim

\demo{Proof} Let $z$ have imaginary part $t$ and real part lying in
$(\frac 12,2]$.  Consider, for $c>\tfrac 12$ 
$$
\frac{1}{2\pi i} \int_{c-i\infty}^{c+i\infty} -\frac{\zeta^{\prime}}{\zeta}(z+w)\frac{x^w}{w^2} dw 
= \sum_{n\le x}\frac{\Lambda(n)}{n^z} \log (x/n),
$$  
 upon integrating term by term using the Dirichlet series expansion of $-\frac{\zeta^{\prime}}{\zeta}(z+w)$. 
 On the other hand, moving the line of integration to the left and calculating residues this 
 equals 
 $$
 -\frac{\zeta^{\prime}}{\zeta}(z) \log x - \Big(\frac{\zeta^{\prime}}{\zeta}(z)\Big)^{\prime} 
 -\sum_{\rho} \frac{x^{\rho-z}}{(\rho-z)^2}+O\Big(\frac{1}{T}\Big). 
 $$
 Integrating from $z=\sigma+it$ to $z=2+it$ we obtain that 
 $$
 \sum_{n\le x} \frac{\Lambda(n)}{n^{\sigma+it} \log n} \log (x/n) + O(1) = 
 (\log x) \log \zeta(\sigma+it) + \frac{\zeta^{\prime}}{\zeta}(\sigma+it) -\sum_{\rho} 
 \int_{\sigma}^{2} \frac{x^{\rho-u-it}}{(\rho-u-it)^2} du.
 $$
The sum over zeros above is bounded in magnitude by 
$$
\sum_\rho \frac{1}{|\rho-\sigma-it|^2} \int_{\sigma}^2 x^{\frac 12-u} du 
\le \frac{x^{\frac 12-\sigma}}{\log x} \sum_{\rho} \frac{1}{|\rho-\sigma-it|^2} 
= \frac{x^{\frac 12-\sigma}}{(\sigma-\frac 12)\log x} F(\sigma+it).
$$
Combining these remarks with (7), the Lemma follows.  
\enddemo 

\proclaim{Lemma 11}  Let $T\le t\le 2T$ be $V$-typical.  
There exists a constant $C$ such that 
for $2\ge \sigma \ge \sigma_0(=\frac 12 +\frac{V}{\log T})$ we have 
$$
\log |\zeta(\sigma+it)| \ge - C V.
$$
\endproclaim
\demo{Proof}  Taking $x=T^{1/V}$ in Lemma 10 and using condition (i) of Definition 3 we 
obtain 
$$
\log |\zeta(\sigma+it)| \ge 
\text{Re } \sum_{n\le x} \frac{\Lambda(n)}{n^{\sigma+it} \log n} \frac{\log (x/n)}{\log x} 
- \frac{2 V }{\log T} F(\sigma+it) 
\ge -2V -\frac{2V}{\log T} F(\sigma+it).
$$
To bound $F(\sigma +it)$, we divide the ordinates $\gamma$ into 
the ranges $2\pi nV/\log T \le |t-\gamma| < 2 \pi (n+1) V/\log T$ for $0\le n\le N=[(\log T)/(4\pi V)]$, 
and the remaining range for $\gamma$. The first kind of zeros contribute, using (ii) of Definition 3, 
$$
\ll V \sum_{0\le n\le N} \frac{(\sigma-\frac 12)}{(\sigma -\frac 12)^2 + (2\pi nV/\log T)^2} \ll \log T.
$$
The remaining zeros contribute 
$$
\sum_{|t-\gamma|>1/2} \frac{(\sigma-\frac12)}{(\sigma-\frac 12)^2 + (t-\gamma)^2} 
\ll \log T.
$$
The Lemma follows. 
\enddemo  

Lemma 11 establishes a stronger form of the desired estimate (3), and as noted 
earlier, the estimate (4) follows from (3) and Lemma 9.  
This completes our proof of Proposition 5.

\enddemo
\head 5. Proof of Theorem 1\endhead 

\noindent We may  assume that $x$ has fractional part half.  
A standard application of Perron's formula (see \S 17 of [1]) gives, with $c= 1+ \frac{1}{\log x}$,
$$
M(x) =\frac{1}{2\pi i} \int_{c-i[x]}^{c+i[x]} \frac{x^s}{s\zeta(s)} ds + O(\log x). \tag{9}
$$
We now deform the contour of integration, replacing the line segment from $c-i[x]$ to 
$c+i[x]$ with a piecewise linear path connecting these points and comprising of a 
number of horizontal and vertical line segments.  We will describe shortly the vertical 
line segments of this contour.  The horizontal line segments shall connect neighboring 
vertical segments, with two end horizontal segments connecting the end vertical 
segments to $c-i[x]$ and $c+i[x]$.  Set $x_0 = [\exp(\sqrt{\log x})]$; one vertical 
segment shall join $\frac 12+ \frac 1{\log x} - ix_0$ to $\frac 12+\frac 1{\log x} +ix_0$.  
For an integer $x_0 \le n \le [x]-1$ we let $V_n$ denote the least integer lying 
in the interval $[(\log \log n)^2, \log n/\log \log n]$ such that all 
points in $[n,n+1]$ are $V_n$-typical.  Notice that the existence of $V_n$ is guaranteed by the 
first assertion of Proposition 4.  
There shall be a vertical 
line segment joining $\frac 12+ \frac{V_n}{\log x} +in$ to $\frac 12 + 
\frac{V_n}{\log x} +i (n+1)$, and its complex conjugate shall also be one of our vertical segments. 
This completes our definition of the contour.

No pole is encountered in deforming our contour, and it remains to estimate the integral 
on these various horizontal and vertical lines.  For the vertical segment from $\frac 12+\frac 1{\log x} -ix_0$ 
to $\frac 12+\frac{1}{\log x}+ix_0$ we use that (see (14.14.2) of [14])
$$ 
|\zeta(\tfrac 12+ \tfrac{1}{\log x} +it)| \gg (|t|+2)^{-\log \log x}
$$
so that 
$$
\Big| \int_{\frac 12+\frac{1}{\log x} -ix_0}^{\frac 12+ \frac{1}{\log x} +ix_0} 
\frac{x^s}{s\zeta(s)} ds \Big| \ll x^{\frac 12} \exp((\log x)^{\frac 12}\log \log x). \tag{10} 
$$

Now suppose $x_0 \le n \le [x]-1$.   The corresponding vertical integral is, 
using Proposition 5,  
$$
\ll \frac{x^{\frac 12}}{n}  \exp(V_n) \exp\Big(V_n \log\frac{\log x}{\log n} + 8 V_n \log\log V_n\Big). \tag{11}
$$ 
Naturally, the same bound applies to the complex conjugate vertical line segment.  Now 
consider the horizontal line segment going from $\frac{1}{2} +\frac{V_n}{\log x}+i(n+1)$ 
to $\frac 12 +\frac{V_{n+1}}{\log x} + i (n+1)$ (if $n=[x]-1$ then the horizontal line segment 
goes from $\frac 12+ \frac{V_n}{\log x} +i[x]$ to $c+i[x]$).  This contributes an amount 
$$
\ll \frac{x^{\frac 12}}{n} \Big( \exp\Big( V_n \log \frac{\log x}{\log n} + 9 V_n\log\log V_n\Big) 
+ \exp\Big( V_{n+1} \log \frac{\log x}{\log (n+1)} +9V_{n+1}\log\log V_{n+1} \Big)\Big). \tag{12}
$$

We split the range for $n$ into dyadic blocks.  Suppose $[T,2T]$ is such 
a dyadic block.  Summing the estimates (11, 12) over elements $n$ in this 
dyadic block we obtain
$$
\ll \frac{x^{\frac 12}}{T} \sum\Sb \log T/\log \log T \ge V \\ V \ge (\log\log T)^2\endSb 
\exp\Big( V\log \frac{\log x}{\log T} +9V\log \log V\Big) \#\{ T\le n\le 2T: \ V_n=V\}. \tag{13}
$$
The terms $V\le 2(\log \log T)^2$ contribute an amount 
$$
\ll x^{\frac 12} \exp\Big( 2(\log \log T)^2 \log \frac{\log x}{\log T} + 18 (\log \log T)^2 \log \log \log T\Big) 
\ll x^{\frac 12} \exp( (\log \log x)^4), 
$$
which is acceptable.  Consider now the contribution of larger values of $V$.  
If $V_n=V$ then by the minimality of $V_n$, it follows that some 
point in $[n,n+1]$ is $(V_n-1)$-atypical.  Appealing to Proposition 4
(pick points from every other interval in order to ensure well-spacing) 
we conclude that the number of such $n$ is $\ll T \exp(- (V-1)\log (V/\log \log T) + 4V\log \log V)$.  
Therefore the quantity in (13) is 
$$
\align
\ll &x^{\frac 12} \exp((\log \log x)^4) \\ 
&+x^{\frac 12} (\log x)^2 \sum\Sb \log T/\log \log T \ge V\\ 
V \ge (\log \log T)^2\endSb  
\exp\Big( V\log \frac{\log x}{\log T} -  V\log \frac{V}{\log \log T} + 13 V\log \log V\Big).
\\
\endalign
$$
A little calculus shows that this is 
$$
\ll x^{\frac 12} \exp((\log \log x)^4) + x^{\frac 12} \exp \Big(\frac{(\log x)\log \log T}{\log T} 
\Big(\log \frac{(\log x)\log \log T}{\log T}\Big)^{13} \Big).
$$
Since $x_0 \le T\le x$, we conclude that the contribution of these 
horizontal and vertical line segments is $\ll x^{\frac 12}\exp((\log x)^{\frac 12} (\log \log x)^{14})$.  
Combining this with (9) and (10) we have established the Theorem.

\enddocument